\documentclass[12pt, reqno]{amsart}
\allowdisplaybreaks
\usepackage{bm,amssymb,amsthm,amsmath}
\usepackage[numbers,sort&compress]{natbib}

\title[]{Global well-posedness of three-dimensional Navier-Stokes equations with partial viscosity under helical symmetry}
\author{Jitao Liu,\,Dongjuan Niu}
\address[Jitao Liu]{College of Applied Sciences, Beijing University of Technology, Beijing, 100124, P. R. China.}
\email{jtliu@bjut.edu.cn,\,\,\,jtliumath@qq.com}
\address[Dongniu Niu]{School of Mathematical Sciences, Capital Normal University, Beijing, 100048, P. R. China.}
\email{djniu@cnu.edu.cn}


\theoremstyle{plain}
\newtheorem{corollary}{Corollary}[section]
\newtheorem{theorem}{Theorem}[section]
\newtheorem{lemma}{Lemma}[section]

\theoremstyle{definition}
\newtheorem{definition}{Definition}[section]

\newtheorem{remark}{Remark}[section]

\let\f=\frac
\let\p=\partial
\def\R{\Bbb R}

\def\no{\noindent}
\def\endproof{\hphantom{MM}\hfill\llap{$\square$}\goodbreak}
\newcommand{\beq}{\begin{equation}}
\newcommand{\eeq}{\end{equation}}
\newcommand{\ben}{\begin{eqnarray}}
\newcommand{\een}{\end{eqnarray}}
\newcommand{\beno}{\begin{eqnarray*}}
\newcommand{\eeno}{\end{eqnarray*}}

\begin{document}


\begin{abstract}
In this paper, we investigate the global well-posedness of three-dimensional Navier-Stokes equations with horizontal viscosity under a special symmetric structure: helical symmetry. More precisely, by a revised Ladyzhenskaya-type inequality and utilizing the behavior of helical flow, we prove the global existence and uniqueness of weak and strong solution to the three-dimensional helical flows. Our result reveals that for the issue of global well-posedness of the viscous helical fluids, the horizontal viscosity plays the important role. To some extent, our work  can be seen as a generalization of the result by Mahalov-Titi-Leibovich [Arch.
Ration. Mech. Anal. 112 (1990), no. 3, 193-222].
\end{abstract}
\maketitle

\section{Introduction and main results}
In present paper, we are concerned with the three-dimensional Navier-Stokes equations  with {\it horizontal viscosity}, which can be read as:
\begin{equation}\label{cauchy}
\left\{\begin{array}{ll}
\p_t\mathbf{u}+\mathbf{u}\cdot\nabla \mathbf{u}-\nu\Delta_{h}\mathbf{u}+\nabla p=0,\\
\nabla\cdot \mathbf{u}=0,
\end{array}\right.
\end{equation}
in a bounded domain $\mathcal{D}\subset\ \mathbb{R}^3$, where $\Delta_{h}=\p^2_{x_1}+\p^2_{x_2}$, $\mathbf{u}=(u_1,u_2,u_3)$ represents the velocity fields, $\nu>0$ is the kinematic
viscosity and $p$ is a scalar pressure. Models with a vanishing anisotropic viscosity in the vertical direction are of relevance for the study of turbulent fows in geophysics. Turbulence is the time-dependent chaotic behavior seen in many fluid flows. This motivates us to study the mathematical problems of fluid flows only with {\it horizontal viscosity}. Starting from Danchin and Paicu \cite{Danchin}, who proved the global well-posedness for the two-dimensional Boussinesq equations only with horizontal viscosity or thermal diffusivity, the topic in this field has attracted considerable attention and great progress has been achieved. One can refer to \cite{CaoWu1,CaoWu3,WZ,JT1,MZ1,MZ2} for details.

As we know, the global well-posedness of three-dimensional Navier-Stokes equations with partial viscosity is far from being resolved. Therefore, we intend to investigate the global well-posedness of solutions to (\ref{cauchy})
under a special symmetric case: helical symmetry. In particular, the flow with helical symmetry is so-called {\it helical flow}, i.e., the flow keeps invariant under certain one-dimensional subgroups of the group of rigid transformations in $\R^3$. These subgroups are generated by a simultaneous rotation around a symmerty axis and a translation along the same symmetry axis. Namely, the subgroup $G^{\kappa}$ is a one-parameter group of isometries of $\R^3$ as follow,
$$G^k = \{S_\theta : \R^3\longrightarrow\R^3 |\theta\in \R\},$$
where $S_\theta$ is the transformation  defined by
\begin{align}\label{trans1}
S_\theta(x)=R_\theta(x)+\left( {\begin{array}{*{20}c}
	0 \\
	0 \\
	\kappa\theta\\
	\end{array} } \right)=\left( {\begin{array}{*{20}c}
	x_1 cos\theta + x_2 sin\theta\\
	-x_1 sin\theta + x_2 cos\theta\\
	x_3+\kappa\theta\\
	\end{array} } \right).	
\end{align}
The nonzero constant $\kappa$  denotes the length scale  and $R_\theta$ is the rotation matrix by an angle $\theta$ around the $x_3$-axis, i.e.,
\begin{align}\label{trans}
R_\theta=\left( {\begin{array}{*{20}c}
	cos\theta& sin\theta&0   \\
	-sin\theta&cos\theta&0 \\
	0&0&1\\
	\end{array} } \right).	
\end{align}
Besides, we will assume $\kappa\equiv1$ for the sake of simplicity throughout the rest paper. As a matter of fact, the
transformation $S_\theta$ corresponds to the superposition of a simultaneous rotation around the $x_3$-axis and a translation along the same $x_3$-axis. The
symmetry lines (orbits of $G^{\kappa}$) are concentric helices. We call the solutions, and more general functions, which are invariant under $G^{\kappa}$ as ``helical''.

Due to the special behavior of helical flow, we intend to discuss our problem in a so-called {\it helical domain},
which is invariant
under the action of $G^{\kappa}$, i.e.,
$$S_\theta\mathcal{D} = \mathcal{D}\quad\quad \forall\theta\in \R.$$
 By recalling that $S_{2\pi}$ is a translation by $2\pi\kappa$ in the $x_3$-direction, we obviously find that
 helical flows inherit a periodic boundary condition in the $x_3$-direction.
 Specially, throughout this paper, we will study system \eqref{cauchy} in a helical domain $\mathcal{D}=\{(x_1,x_2,x_3)\in\R^3\,|\,0<x_1^2+x_2^2<1,0<x_3<2\pi\}=B_1\times[0,2\pi]$ and implement the following initial-boundary condition to \eqref{cauchy}
\begin{equation}\label{cauchy0}
\left\{\begin{array}{ll}
\mathbf{u}(x_1,\,x_2,\,x_3,\,t)=\mathbf{u}(x_1,\,x_2,\,x_3+2\pi,\,t),\ \ \,x\in\mathcal{D},\\
\mathbf{u}(x_1,\,x_2,\,x_3,\,t)|_{\p\mathcal{D}}=0,\\
\mathbf{u}(x,0)=\mathbf{u}_0(x).\\
\end{array}\right.
\end{equation}

As common in practice as axisymmetric flows, helical flows have attached wide mathematical attention recently. In 1990, global existence and uniqueness of strong solutions to three-dimensional Navier-Stokes equations with helical symmetry have been obtained by Mahalov, Titi and Leibovich \cite{MTL} with helical initial data.
 For Euler equations with helical symmetry, when the helical swirl (which is similar to the quantity in the axisymmetric case) vanished, Dutrifoy \cite{Du}  obtained the global existence and uniqueness of smooth solutions in bounded domain with the regular initial
data. Ettinger and Titi \cite{ET} derived the global existence and uniqueness
of strong solutions in bounded domain for bounded initial vorticity. In addition, for weaker initial assumptions, the existence of weak solutions to Euler equations with helical symmetry  has also been discussed in \cite{BLL} and \cite{JLN}.
In \cite{BMNHT}, the authors proved the stability of weak solutions to the three-dimensional Navier-Stokes with helical initial data. In \cite{MLLNT}, when the helical parameter $\kappa$ goes to infinity, the limiting property of the incompressible flows with helical symmetry has been investigated.

In this paper, we are devoted to studying the global existence, uniqueness and stability of weak and strong solutions to \eqref{cauchy}-\eqref{cauchy0} with helical symmetry. In our case, on account of vertical smoothing effect vanishing in the equations \eqref{cauchy}, the viscosity term is not enough to control the nonlinear term in the energy estimates. To overcome this difficulty, inspired by the key observation of \cite{MTL}, we establish a {\it revised} version of Ladyzhenskaya-type inequality,
\ben\label{key02}
\|\mathbf{u}\|_{L^4(\mathcal{D})}\leq \sqrt[4]{\f2{\pi}}\big[\|\mathbf{u}\|_{L^2(\mathcal{D})}^{\f12}\|\mathbf{\nabla_{h}u}\|_{L^2(\mathcal{D})}^{\f12}+\|\mathbf{u}\|_{L^2(\mathcal{D})}\big],
\een
which gives us the cornerstone to obtain the global existence of weak solution. However, to prove the stability of weak solution and global existence of strong solution,
we still have some difficulties in deriving the necessary estimate of $\|\nabla \mathbf{u}\|_{L^2(0,t;L^2({\mathcal{D}}))}$ (independent of time $t$) even with the help of (\ref{key02}). To get over this problem, we make full advantage of the helically symmetric structure of the flows. First, we establish the exponential decay estimates of velocity field. Then, thanks to a novel observation (see \eqref{novel} below), we set up the desired estimates $\|\nabla \mathbf{u}\|_{L^2(0,t;L^2({\mathcal{D}}))}$, which is uniform with respect to the time variable. At the end, we derive the estimate of
$\|\mathbf{u}\|_{L^\infty(0,t;H^1({\mathcal{D}}))}$, which helps us to resolve the problem.

Before showing the main theorems, we would like to introduce the spaces $H(\mathcal{D})$ and $V(\mathcal{D})$ be the closure of $C^{\infty}$ vector fields which are periodic in the vertical variable, compactly supported in the horizontal sections, and divergence-free in $\mathcal{D}$ with respect to $L^2(\mathcal{D})$ and $H^1_{\rm ver}(\mathcal{D})=\{\mathbf{u}\in L^2(\mathcal{D}),\,\,\nabla_h \mathbf{u}\in L^2(\mathcal{D})\}$ norms, respectively. On this basis, we then define the inner products of $L^2(\mathcal{D})$ by $(\mathbf{u},\mathbf{v})=\sum\limits_{{i=1}}^3\int_{\mathcal{D}}u_iv_i\,dx$, and   denote by $V'$ the dual space of $V$ and the action of $V'$ on $V$ by $<\cdot\,,\,\cdot>$. Moreover, we use the following notation for the trilinear continuous form by setting
\ben\label{b0}
b(\mathbf{u},\mathbf{v},\mathbf{w})=\sum\limits_{{i,j=1}}^3\int_{\mathcal{D}}u_i\p_iv_jw_j\,dx.
\een
If $u\in\,V$, then
\ben\label{b1}
b(\mathbf{u},\mathbf{v},\mathbf{w})=-b(\mathbf{u},\mathbf{w},\mathbf{v}),\,\,\,\forall\,\,v,w\in V(\mathcal{D}),
\een
and
\ben\label{b2}
b(\mathbf{u},\mathbf{v},\mathbf{v})=0,\,\,\,\forall\,\,v\in V(\mathcal{D}).
\een
Finally, we will set up the definitions of weak and strong solutions of system (\ref{cauchy})-(\ref{cauchy0}) as bellow.

\begin{definition}\label{weak}
({\bf Weak solution}):\,
Suppose $\mathbf{u_0}\in H(\mathcal{D})$ be helically symmetric, the helical vector fields $\mathbf{u}(x,t)$ is called a global weak solution of (\ref{cauchy})-(\ref{cauchy0}) if for any $t>0$,
\ben\label{weak1}
\|\mathbf{u}(t,\cdot)\|_{L^2(\mathcal{D})}^2+2\nu\int_0^t\|\nabla_h\mathbf{u}(\tau,\cdot)\|_{L^2(\mathcal{D})}^2d\tau\leq \|\mathbf{u_0}\|_{L^2(\mathcal{D})}^2,
\een
and
\ben\label{weak2}
\int_{\mathcal{D}}\mathbf{u_0}\cdot\pmb{\varphi_0}\,dx+\int_0^t\int_{\mathcal{D}}\big[\mathbf{u}\cdot\pmb{\varphi_t}+\mathbf{u}\cdot\nabla \pmb{\varphi}\cdot \mathbf{u}-\nu\nabla_h \mathbf{u}:\nabla_h\pmb{\varphi}\big]dxd\tau=0,
\een
holds for any helical vector fields $\pmb{\varphi}\in C_{c}^\infty([0,t)\times{\mathcal{D}})$ with $\nabla\cdot\pmb{\varphi}=0$, where $A:B\equiv\sum\limits_{i,j}a_{ij}b_{ij}$ is the trace product of two matrices.
\end{definition}

\begin{remark}\label{gw1}
	Following standard arguments as in the theory of the Navier-Stokes equations (see e.g., \cite{RT}),we mention that system (\ref{weak2}) is equivalent to the following system
	\beno
	&&\f{d}{dt}<\pmb{u},\pmb{\varphi}>+\,b(\mathbf{u},\mathbf{u},\pmb{\varphi})+(\nabla_h \mathbf{u},\nabla_h \pmb{\varphi})=0,\quad\quad\quad\forall\,\pmb{\varphi}\in L^2(0,t;V(\mathcal{D})).
	\eeno
\end{remark}

\begin{definition}\label{strong}
({\bf Strong solution}):\, Let Definition \ref{weak} be satisfied. Furthermore, it holds that
\ben\label{weak3}
\mathbf{u}\in L^\infty(0,t;H^1(\mathcal{D})),\,\,\nabla_{h}\mathbf{u}\in L^2(0,t;H^1(\mathcal{D})),
\een
for any $t>0$, then we call the corresponding solution as a global strong solution.
\end{definition}

Now, we are in the position to state the main results of this paper.

\begin{theorem}\label{gw}
Suppose that $\mathbf{u_0}\in H(\mathcal{D})$, then there exists a unique global weak solution $\mathbf{u}\in L^{\infty}(0,\infty;H(\mathcal{D}))\cap L^{2}(0,\infty;V(\mathcal{D}))$
to (\ref{cauchy})-(\ref{cauchy0}). Moreover,  $\mathbf{u}$ satisfies the exponential decay rate
\ben\label{decay}
\|\mathbf{u}(\cdot,t)\|_{L^2(\mathcal{D})}^2\leq e^{-\f{2\nu}{c_0}t}\|{\mathbf{u_0}}\|_{L^2(\mathcal{D})}^2 \quad \forall  t\geq 0,
\een
where $c_0$ is the constant of Poincar{\'e} inequality.
\end{theorem}


Motivated by the work of Bardos et al. in \cite{BMNHT}, we also observe a corresponding stability result.

\begin{theorem}\label{gsta}
Given that $\mathbf{u_0}\in H(\mathcal{D})$ a helical vector field, there exists 
$$\mathbf{u}\in L^{\infty}(0,\infty;H(\mathcal{D}))\cap L^{2}(0,\infty;V(\mathcal{D})),$$
 the weak solution of helical incompressible Navier-Stokes equations  (\ref{cauchy})-(\ref{cauchy0}) with initial data $\mathbf{u_0}$, given in Theorem
\ref{gw}. Moreover, let $\mathbf{v_0}\in H(\mathcal{D})$ be a general vector field and
$$\mathbf{v}\in L^{\infty}(0,\infty;H(\mathcal{D}))\cap L^{2}(0,\infty;V(\mathcal{D})),$$
be a Leray-Hopf weak solution of three-dimensional incompressible Navier-Stokes equations (\ref{cauchy})-(\ref{cauchy0}) with initial data $\mathbf{v_0}$. Then it holds that
$$\|\mathbf{u}-\mathbf{v}\|_{L^2(\mathcal{D})}^2(t)\leq \|\mathbf{u_0}-\mathbf{v_0}\|_{L^2(\mathcal{D})}^2\,{\rm exp} \big[\f{4c_*^2+2c_*^2c_0}{{\nu^2}}\|\mathbf{u_0}\|_{L^2(\mathcal{D})}^2\big],$$
for any $t\geq0$, where $c_0$ and $c_*$ are the constants from Poincar{\'e} inequality and Lemma \ref{twointegral}, respectively.
\end{theorem}

Moreover, for the strong solution of \eqref{cauchy}-\eqref{cauchy0}, we also have the similar argument as follows.
\begin{theorem}\label{gs}
Assume that $\mathbf{u_0}\in V(\mathcal{D})$, then there exists a unique global strong solution $\mathbf{u}\in L^\infty(0,t;H^1(\mathcal{D}))$ to the system (\ref{cauchy})-(\ref{cauchy0}) in the sense of Definition \ref{strong}.
\end{theorem}

This paper is organized as follows. In section 2, we introduce some notations and technical lemmas. Section 3 is devoted to the a priori estimates and proof of main theorems.

\section{Preliminary}\hspace*{\parindent}

In this section, we will fix some notations and set down some basic definitions. As discussed in the Introduction, helical flow is invariant under certain one-dimensional subgroups. Especiallly, we will employ the following definition.
\begin{definition}(\it Helical flow)\label{def1}\\
	(i):\,\,\,A scalar function $f: \mathbb{R}^3\longrightarrow \mathbb{R}$ is said to be helical if
	\beno
	f(S_\theta(x))=f(x),\quad \forall\,\theta\in\R.
	\eeno
	(ii):\,\,\,A vector field $\mathbf{v}: \mathbb{R}^3\longrightarrow \mathbb{R}^3$ is said to be helical if
	\beno
	\mathbf{v}(S_\theta(x))=R_\theta \mathbf{v}(x),\quad \forall\,\theta\in\R.
	\eeno
\end{definition}
Subsequently, we would like to introduce some important properties of helical flow. By setting $\pmb{\xi}=(x_2,\,-x_1,\,1)^{\top}$, the following lemmas hold.

\begin{lemma}(Claim 2.5, \cite{ET})\label{key1} A smooth vector field $\mathbf{v}=(v_1,\,v_2,\,v_3)^{\top}: \mathbb{R}^3\longrightarrow \mathbb{R}^3$ is helical if
and only if it the following relation holds true:
$$\p_{\pmb{\xi}}\mathbf{v}=\mathbf{v}^{\perp},$$
where $\mathbf{v}^{\perp}=(v_2,\,-v_1,\,0)^{\top}.$
\end{lemma}

The following lemma tells that the helical flow can be essentially viewed as an extension of the two-dimensional one in some sense (see Proposition 2.1 of \cite{MLLNT}), which make it possible to improve the result of \cite{BLL} to the two-dimensional case \cite{Majda}.

\begin{lemma} \label{twodimpro}
Let $\mathbf{u} = \mathbf{u}(x) $ be a smooth helical vector field and
let $p=p(x)$ be a smooth helical function, where
$x=(x_1,x_2,x_3)$. Then there exist {\em unique\/} $\mathbf{w}=(w^1,w^2,w^3)
=(w^1,w^2,w^3)(y_1,y_2)$ and $q=q(y_1,y_2)$ such that
\begin{equation} \label{uandpinY}
\mathbf{u}(x) = R_{2\pi x_3}\mathbf{w}(y(x)), \;\;\; p=p(x)=q(y(x)),\end{equation}
with $R_\theta$ given in \eqref{trans},
and
\begin{equation} \label{Yofx}
y(x) =
\left[\begin{array}{l}
y_1 \\ \\
y_2
\end{array}\right]
=
\left[\begin{array}{cc}
\cos(2\pi x_3/\sigma) & -\sin(2\pi x_3/\sigma) \\ \\
\sin(2\pi x_3/\sigma) & \cos(2\pi x_3/\sigma)
\end{array}\right]
\left[\begin{array}{l}
x_1\\ \\
x_2
\end{array}\right].
\end{equation}

Conversely, if $\mathbf{u}$ and $p$ are defined through \eqref{uandpinY} for some
$\mathbf{w}=\mathbf{w}(y_1,y_2)$, $q=q(y_1,y_2)$,
then $\mathbf{u}$ is a helical vector field and $p$ is a helical scalar function.
\end{lemma}

In the end, we will provide an {\it revised version} of Ladyzhenskaya-type inequality, which is given in \cite{MTL}. Our version would be slightly different and plays an important role in the present paper.
\begin{lemma}(Lemma 3.1, \cite{MTL})\label{key2}
Let ${\mathbf{u}}$ be a helical function in $H^1(\mathcal{D})$, then it follows that
\ben\label{key21}
\|\mathbf{u}\|_{L^4(\mathcal{D})}\leq \sqrt[4]{\f2{\pi}}\,\,\big[\|\mathbf{u}\|_{L^2(\mathcal{D})}^{\f12}\|\nabla_{h}\mathbf{u}\|_{L^2(\mathcal{D})}^{\f12}\big],
\een
where $\nabla_{h}=(\p_{x_1},\,\p_{x_2}).$ if in addition, $\mathbf{u}|_{\p\mathcal{D}}=0$, (\ref{key21}) will be replaced by
\ben\label{key22}
\|\mathbf{u}\|_{L^4(\mathcal{D})}\leq C\,\,\|\mathbf{u}\|_{L^2(\mathcal{D})}^{\f12}\|\nabla_{h}\mathbf{u}\|_{L^2(\mathcal{D})}^{\f12},
\een
where $C$ is the generic constant.
\end{lemma}
\no{\bf Proof.}\quad
According to Lemma \ref{twodimpro},  there exists the corresponding vector field $\boldsymbol{w}$ by virtue of
\eqref{uandpinY}. It is easily to compute that
$|\mathbf{u}(x_1,x_2,x_3)|^2=|\boldsymbol{w}(y_1,y_2)|^2.$
Then
\begin{align}\label{ineq3}
\|\mathbf{u}\|_{L^4(\mathcal{D})}^4
=2\pi\|\boldsymbol{w}\|_{L^4(B_1)}^4
\leq C\|\boldsymbol{w}\|_{L^2(B_1)}^2\|\nabla_y \boldsymbol{w}\|_{L^2(B_1)}^2.
\end{align}
In the above inequality,  we use the two-dimensional Sobolev embedding inequality.
Again using the direct equation \eqref{uandpinY},  we deduce that
\begin{align}\label{ineq1}
\|\boldsymbol{w}\|_{L^2(B_1)}= \frac{1}{\sqrt{2\pi}}\|\mathbf{u}\|_{L^2(\mathcal{D})},
\end{align}
and
\begin{align}\label{ineq2}
\|\nabla_y \boldsymbol{w}\|_{L^2(B_1)} \leq \|\nabla_h \mathbf{u}\|_{L^2(\mathcal{D})}.
\end{align}
Then,  taking \eqref{ineq1} and \eqref{ineq2} into \eqref{ineq3}, we finish the proof.
\endproof\vskip 0.5cm

\begin{lemma} \label{twointegral}
Let $\mathbf{u} = \mathbf{u}(x)\in H^1(\mathcal{D}) $ be a smooth helical vector field ,then there exist a constant $c_{*}$ such that
$$\|\nabla \mathbf{u}\|_{L^2(B_1)}\leq c_{*}\| \mathbf{u}\|_{H^1(\mathcal{D})}.$$
\end{lemma}
\no{\bf Proof.}\quad
As $u$ is helical flow, it is evident that the $L^2({\mathcal{D}})$ norm of $\nabla_h u$ is independent of $x_3$. Specially,  based on the two-dimensional property of Lemma \ref{twodimpro}, there exists the corresponding vector field
$\boldsymbol{w} \in L^{\infty}(0,T;H^1_{loc}(\mathbb{R}^2)$ by virtue of \eqref{uandpinY}.  Then direct computations indicate that each component of $\nabla_y \boldsymbol{w}$ is a composition of the components of $\partial_{x_1} \mathbf{u}, \partial_{x_2} \mathbf{u}$
and the trigonometric functions about $x_3$ according to the formulae of \eqref{uandpinY} and \eqref{Yofx}.

Without loss of generality, we take the expression of  $\partial_{y_1} w_1$ for instance, \textit{i.e.},
$$\partial_{y_1}w_1=\partial_{x_1} u_1\cos^2 x_3-\partial_{x_1} u_2\cos x_3\sin x_3-\partial_{x_2} u_1\cos x_3\sin x_3+\partial_{x_2}v_2\sin^2 x_3.$$
Therefore, it is easy to deduce that for any $t\in (0,T)$
$$c\|\nabla_h \mathbf{u}\|_{L^2(B_1)}\leq \|\nabla_y \boldsymbol{w}\|_{L^2(B_1)}\leq C\|\nabla_h \mathbf{u}\|_{L^2(\mathcal{D})},$$
and
$$
c\|\mathbf{u}\|_{L^2(B_1)}\leq \|\boldsymbol{w}\|_{L^2(B_1)}\leq C\|\mathbf{u}\|_{L^2(\mathcal{D})},
$$
where $c$ and $C$ are  generic constants.
Moreover,  thanks to the helical property of Lemma \ref{key1},
$$\|\partial_{x_3} \mathbf{u}\|_{L^2(B_1)}\leq \|\nabla_h \mathbf{u}\|_{L^2(B_1)}+\|\mathbf{u}\|_{L^2(B_1)}.$$

As a result, one can find a constant $c_{*}$ such that for any $t\in (0,T)$
$$\|\nabla \mathbf{u}\|_{L^2(B_1)}\leq c_{*}\| \mathbf{u}\|_{H^1(\mathcal{D})}.$$
\endproof\vskip 0.5cm
\section{Global well-posedness}\hspace*{\parindent}

\subsection{{\em A priori} estimates}\hspace*{\parindent}\\

In this subsection, we will establish the {\em a priori} estimates of velocity fields. At first, we will list the basic energy estimate with decay rate.
\begin{lemma}\label{lem1}
Suppose $\mathbf{u_0}\in {L^2({\mathcal{D}})}$ be a helical function with $\nabla\cdot \mathbf{u_0}=0$, then for a helical smooth solution $\mathbf{u}$ of (\ref{cauchy}) and any $t\geq0$, there holds that
\ben\label{energy1}
\|{\mathbf{u}}(\cdot,t)\|_{L^2(\mathcal{D})}^2\leq e^{-\f{2\nu}{c_0}t}\|{\mathbf{u_0}}\|_{L^2(\mathcal{D})}^2,
\een
and
\ben\label{energy12}
\nu\int_0^te^{\f{\nu}{c_0}\tau}\|\nabla_h\mathbf{u}(\tau)\|^{2}_{L^2(D)} d\tau\leq\|\mathbf{u_0}\|^{2}_{L^2(\mathcal{D})},
\een
where $c_0$ is the constant in Poincar{\'e} inequality.
\end{lemma}
\no{\bf Proof.}\quad
Taking inner product of $(\ref{cauchy})^1$ with $\mathbf{u}$, and then integrating over ${\mathcal{D}}$, it follows that
\ben\label{lem11}
\f{d}{dt}\|\mathbf{u}\|_{L^2(\mathcal{D})}^2+2\nu\|\nabla_{h}\mathbf{u}\|_{L^2(\mathcal{D})}^2\leq\ 0.
\een
Moreover, by noticing that $\mathbf{u}|_{\p\mathcal{D}}=0$ and implying Poincar{\'e} inequality, we have
\ben\label{lem12}
\|\mathbf{u}\|_{L^2(\mathcal{D})}^2\leq c_0\|\nabla_{h}\mathbf{u}\|_{L^2(\mathcal{D})}^2
\een
for some constant $c_0$. Then, we can rewrite \eqref{lem11} as
\beno
\f{d}{dt}\|\mathbf{u}\|_{L^2(\mathcal{D})}^2+\f{2\nu}{c_0}\|\mathbf{u}\|_{L^2(\mathcal{D})}^2\leq\ 0,
\eeno
which yields
\ben\label{lem13}
\|\mathbf{u}(\cdot,t)\|_{L^2(\mathcal{D})}^2\leq e^{-\f{2\nu}{c_0}t}\|\mathbf{u_0}\|_{L^2(\mathcal{D})}^2 ,\quad\forall t\geq0.
\een

Then, by multiplying \eqref{lem12} with $e^{\f{\nu}{c_0}t}$ and usig \eqref{lem13}, we have
\begin{equation*}\begin{split}
&\frac12\frac{d}{dt}e^{\f{\nu}{c_0}t}\|\mathbf{u}\|^{2}_{L^2(\mathcal{D})}+\nu e^{\f{\nu}{c_0}t} \|\nabla_h\mathbf{u}\|^{2}_{L^2(\mathcal{D})}\\
&\leq\f{\nu}{2c_0}e^{\f{\nu}{c_0}t}\|\mathbf{u}\|^{2}_{L^2(\mathcal{D})}\\
&\leq\f{\nu}{2c_0}e^{-\f{\nu}{c_0}t}\|\mathbf{u_0}\|^{2}_{L^2(\mathcal{D})},
\end{split}
\end{equation*}
which also yields, after integrating in time over $[0,t]$, that
\begin{equation}\label{lem15}
\begin{split}
\nu\int_0^te^{\f{\nu}{c_0}\tau}\|\nabla_h\mathbf{u}(\tau)\|^{2}_{L^2(D)} d\tau\leq\|\mathbf{u_0}\|^{2}_{L^2(\mathcal{D})}.
\end{split}
\end{equation}
\endproof\vskip 0.5cm

Nex, we make an attempt to derive the $H^1$ estimate of $\mathbf{u}$, especially the uniform bound of $\|\nabla \mathbf{u}\|_{L^2(0,t;L^2({\mathcal{D}}))}$. Noicing that there is no vertical smoothing effect in the system \eqref{cauchy}, we only have the estimate of $\|\nabla_h \mathbf{u}\|_{L^2(0,t;L^2({\mathcal{D}}))}$ other than $\|\nabla \mathbf{u}\|_{L^2(0,t;L^2({\mathcal{D}}))}$ in \eqref{energy1}. This gives rise to the difficulty in deriving this estimate of $\|\mathbf{u}\|_{L^2(0,t;H^1({\mathcal{D}}))}$, even with the help of {\it revising } Ladyzhenskaya-type inequality. To overcome it, we make full advantage of the helical symmetry structure, which yields the following conclusion.

\begin{corollary}\label{cor1}
Let $\mathbf{u_0}$ be as in Lemma \ref{lem1}, then for a smooth helical solution $\mathbf{u}$ of (\ref{cauchy}) and any $t\geq0$, there holds
\ben\label{energy2}
\int_0^t\|\nabla \mathbf{u}\|_{{L^2({\mathcal{D}})}}^2d\tau\leq\f{4+c_0}{2\nu}\|{\mathbf{u_0}}\|_{L^2(\mathcal{D})}^2,
\een
where $c_0$ is the constant in Poincar{\'e} inequality.
\end{corollary}
\no{\bf Proof.}\quad
Thanks to Lemma \ref{key1}, for any helical flow, it holds that
\ben\label{novel}
\p_{x_3}\mathbf{u}={x_2}\p_{x_1}\mathbf{u}-{x_1}\p_{x_2}\mathbf{u}+\mathbf{u}^{\bot}.
\een
Then, we notice that $\mathcal{D}$ is bounded by 1 in $x_1$ and $x_2$ direction. Thus, by applying \eqref{lem13} and \eqref{lem15}, we have
\beno
\int_0^t\|\nabla \mathbf{u}\|_{{L^2({\mathcal{D}})}}^2d\tau&\leq&\int_0^t\|\p_{x_3} \mathbf{u}\|_{{L^2({\mathcal{D}})}}^2d\tau+\int_0^t\|\nabla_{h}\mathbf{u}\|_{{L^2({\mathcal{D}})}}^2d\tau\\
&\leq&\int_0^t\|\mathbf{u}\|_{{L^2({\mathcal{D}})}}^2d\tau+2\int_0^t\|\nabla_{h}\mathbf{u}\|_{{L^2({\mathcal{D}})}}^2d\tau\\
&\leq&[\int_0^te^{-\f{2\nu}{c_0}\tau}d\tau+2{\nu}^{-1}]\|{\mathbf{u_0}}\|_{L^2(\mathcal{D})}^2\\
&\leq&\f{4+c_0}{2\nu}\|{\mathbf{u_0}}\|_{L^2(\mathcal{D})}^2.
\eeno
\endproof

On the basis of Corollary \ref{cor1}, we have derived the  estimate $\|\nabla \mathbf{u}\|_{L^2(0,t;L^2({\mathcal{D}}))}$ independent  of $t$. This gives us a good cornerstone to estimate $\|\mathbf{u}\|_{L^\infty(0,t;H^1({\mathcal{D}}))}$, which is necessary to prove the existence of strong solutions (Theorem \ref{gs}).

\begin{lemma}\label{lem3}
Under the assumptions of Theorem \ref{gs}, for a smooth helical solution $\mathbf{u}$ of (\ref{cauchy}), there holds
\ben\label{estimate31}
\|\mathbf{u}\|_{H^1(\mathcal{D})}^2+\nu\int_0^t\|\nabla_{h}\mathbf{u}\|_{H^1(\mathcal{D})}^2d\tau\leq\ C,
\een
where the constant $C$ depends only on t.
\end{lemma}
\no{\bf Proof.}\quad
Taking inner product of $(\ref{cauchy})^1$ with $-\Delta \mathbf{u}$ on $\mathcal{D}$, integrating by parts, employing the boundary condition (\ref{cauchy0}) and applying H\"{o}lder inequalities, Lemma \ref{key2}, we have
\beno
&&\f12\f{d}{dt}\|\nabla \mathbf{u}\|_{L^2({\mathcal{D}})}^2+\nu\|\Delta_h \mathbf{u}\|_{L^2({\mathcal{D}})}^2+\nu\|\p_{x_3}\nabla_h \mathbf{u}\|_{L^2({\mathcal{D}})}^2\\
&=&-\int_{\mathcal{D}}\mathbf{u}\cdot\nabla \mathbf{u}\cdot \Delta \mathbf{u}dx\\
&=&\int_{\mathcal{D}}\sum_{i=x,y,z}\mathbf{u}\cdot\nabla \p_i\mathbf{u}\cdot \p_i \mathbf{u}dx+\int_{\mathcal{D}}\sum_{i=x,y,z}\p_i\mathbf{u}\cdot\nabla \mathbf{u}\cdot \p_i \mathbf{u}dx\\
&=&\int_{\mathcal{D}}\sum_{i=x,y,z}\p_i\mathbf{u}\cdot\nabla \mathbf{u}\cdot \p_i \mathbf{u}dx\\
&\leq&C\|\nabla \mathbf{u}\|_{L^4({\mathcal{D}})}^2\|\nabla \mathbf{u}\|_{L^2({\mathcal{D}})}\\
&\leq&C\|\nabla \mathbf{u}\|_{L^2({\mathcal{D}})}^2(\|\nabla_h^2 \mathbf{u}\|_{L^2({\mathcal{D}})}+\|\p_{x_3} \nabla_h \mathbf{u}\|_{L^2({\mathcal{D}})})+C\|\nabla \mathbf{u}\|_{L^2({\mathcal{D}})}^3.
\eeno
Then by using the elliptic theory and boundary condition (\ref{cauchy0}) again, one can derive that $\|\nabla_h^2 \mathbf{u}\|_{L^2({\mathcal{D}})}\leq C\|\Delta_h \mathbf{u}\|_{L^2({\mathcal{D}})}$. Thanks to it and Young inequalities, one has
\beno
&&\f12\f{d}{dt}\|\nabla \mathbf{u}\|_{L^2({\mathcal{D}})}^2+\|\Delta_h \mathbf{u}\|_{L^2({\mathcal{D}})}^2+\|\p_{x_3}\nabla_h \mathbf{u}\|_{L^2({\mathcal{D}})}^2\\
&\leq&C\|\nabla \mathbf{u}\|_{L^2({\mathcal{D}})}^2(\|\Delta_h \mathbf{u}\|_{L^2({\mathcal{D}})}+\|\p_{x_3} \nabla_h \mathbf{u}\|_{L^2({\mathcal{D}})})+C\|\nabla \mathbf{u}\|_{L^2({\mathcal{D}})}^3.\\
&\leq&\f{\nu}2[\|\Delta_h \mathbf{u}\|_{L^2({\mathcal{D}})}^2+\|\p_{x_3}\nabla_h \mathbf{u}\|_{L^2({\mathcal{D}})}^2]+C[1+\|\nabla \mathbf{u}\|_{L^2({\mathcal{D}})}^4],
\eeno
which implies the conclusion after applying Gronwall's inequality and Corollary \ref{cor1}.
\endproof\vskip 0.5cm

\subsection{Proof of Theorem \ref{gw}}\hspace*{\parindent}\vskip 0.5cm

The proof will be devided into three steps, that is, {\it existence and uniqueness}.\vskip 0.5cm

\no{\bf Existence:} It is well-known that there exists at least a weak solution
of (\ref{cauchy}) provided that the initial velocity belongs to $H(\mathcal{D})$. The proof
is standard (see e.g. [17]), and we only list a sketch here. For $\epsilon>0$, let $j$ be a positive radial compactly supported smooth function whose integral equals 1 and denote $J_\epsilon$ as a Friedrichs mollifier by
\ben\label{mollifier}
J_\epsilon \mathbf{u}=j_\epsilon\ast \mathbf{u},\quad\,\hbox{where}\,\,\,\,\,j_\epsilon=\epsilon^{-2}j(\epsilon^{-1}x).
\een
Moreover, let $\mathcal{P}$ be the Leray projector over divergence free vector fields, then the following properties holds
\ben\label{jp}
J^2_\epsilon=J_\epsilon,\quad\,\mathcal{P}^2=\mathcal{P},\quad\,\mathcal{P}J_\epsilon=J_\epsilon\mathcal{P}.
\een

Then, we construct the following approximating system
\begin{equation}\label{cauchy3}
\left\{\begin{array}{ll}
\p_t\mathbf{u}_\epsilon+\mathcal{P}J_\epsilon(J_\epsilon \mathbf{u}_\epsilon\cdot\nabla J_\epsilon \mathbf{u}_\epsilon)-\nu\Delta_{h}\mathcal{P}J_\epsilon \mathbf{u}_\epsilon=0,\\
\nabla\cdot \mathbf{u}_\epsilon=0,\\
\mathbf{u}_\epsilon(x,0)=J_\epsilon \mathbf{u}_0.
\end{array}\right.
\end{equation}
By the Picard theorem (see e.g. \cite{Majda}), there exists a unique solution of (\ref{cauchy}). Similar to the standard method, we can prove that
$\{\mathbf{u}_{\epsilon}\}$ is a Cauchy sequence for any $\epsilon>0$. Futhermore, by Proposition 1.1 in \cite{Majda}, for Euler and Navier-Stokes equations, the transformations {\it Galilean invariance, Rotation symmetry, Scale invariance} also yields solutions.  This shows that $\{\mathbf{u}_{\epsilon}\}$ also preserves the helical symmetry. Finally, it is not hard to prove that the limit $\mathbf{u}$ of sequence $\{\mathbf{u}_{\epsilon}\}$ satisfies (\ref{cauchy}) in the distribution sense as $\epsilon$ tends to zero.\vskip 0.5cm

\no{\bf Uniqueness:} For any fixed $t>0$, suppose there are two solutions $(\mathbf{u},p)$, $(\mathbf{\widetilde{u}},\widetilde{p})$ of (\ref{cauchy}) and let $\mathbf{U}=\mathbf{\widetilde{u}}-\mathbf{u},\,\,P=\widetilde{p}-p,$ then by Remark \ref{gw1}, it holds that
\ben\label{unix}
\f{d}{dt}<\mathbf{U},\pmb{\varphi}>+\,b(\mathbf{\widetilde{u}},\mathbf{U},\pmb{\varphi})+b(\mathbf{U},\mathbf{u},\pmb{\varphi})+(\nabla_h\mathbf{U},\nabla_h\pmb{\varphi})=0,
\een
\ben\label{uniz}
\mathbf{U}(x,0)=0,
\een
for any $\pmb{\varphi}\in L^2(0,t;V(\mathcal{D}))$.

Subsequently, by taking $\pmb{\varphi}=\mathbf{U}$ in $(\ref{unix})$, making use of (\ref{b2}) and Lions-Magenes Lemma (see e.g., \cite{RT}), we have
\begin{align}
\f12\f{d}{dt}\|\mathbf{U}\|_{L^2(\mathcal{D})}^2+\|\nabla_h \mathbf{U}\|_{L^2(\mathcal{D})}^2&
\leq-\int_{\mathcal{D}}\mathbf{U}\cdot\nabla \mathbf{u}\cdot \mathbf{U}\,d\nonumber x\\
&\leq C\|\mathbf{U}\|_{L^4(\mathcal{D})}^2\|\nabla \mathbf{u}\|_{L^2(\mathcal{D})}\nonumber\\
&\leq C\|\mathbf{U}\|_{L^2(\mathcal{D})}\|\nabla_h \mathbf{U}\|_{L^2(\mathcal{D})}\|\nabla \mathbf{u}\|_{L^2(\mathcal{D})}\\
&\leq\f12\|\nabla_h \mathbf{U}\|_{L^2(\mathcal{D})}^2+C\|\nabla \mathbf{u}\|_{L^2(\mathcal{D})}^2\|\mathbf{U}\|_{L^2(\mathcal{D})}^2.\nonumber
\end{align}

In the end, by employing Corollary \ref{cor1} and Gronwall's inequality, it yields that
\beno
\|\mathbf{U}(\tau)\|_{L^2(\mathcal{D})}^2\leq C\|\mathbf{U}_0\|_{L^2(\mathcal{D})}^2=0,
\eeno
for any $\tau\in[0,t].$ Thus, the proof is finished.
\endproof\vskip 0.5cm

\subsection{Proof of Theorem \ref{gsta}}\hspace*{\parindent}\vskip 0.5cm

To begin with, for any fixed $t>0$ and $j_{\epsilon}$ (defined in (\ref{mollifier})), we define the notation $f^\epsilon$ as below
$$f^\epsilon=f^\epsilon(t,x)=\int_0^tj_{\epsilon}(t-\tau)f(\tau,x)d\tau.$$
Then we choose $\mathbf{v}^{\epsilon}$ and $\mathbf{u}^{\epsilon}$ as the test function in the weak formulation for $\mathbf{u}$ and $\mathbf{v}$ respectively. Furthermore,  it follows that
\ben\label{sta1}
&&-(\mathbf{u},J_{\epsilon}\mathbf{v})(t)+\int_0^t(\mathbf{u},\p_tJ_{\epsilon}\mathbf{v})d\tau-\nu\int_0^t(\nabla_h \mathbf{u},\nabla_h J_{\epsilon}\mathbf{v})d\tau\\
&&=-\int_0^t(\mathbf{u}\cdot\nabla J_{\epsilon}\mathbf{v}, \mathbf{u})d\tau-(\mathbf{u_0},J_{\epsilon}\mathbf{v_0})\notag,
\een
and
\ben\label{sta2}
&&-(\mathbf{v},J_{\epsilon}\mathbf{u})(t)+\int_0^t(\mathbf{v},\p_tJ_{\epsilon}\mathbf{u})d\tau-\nu\int_0^t(\nabla_h \mathbf{v},\nabla_h J_{\epsilon}\mathbf{u})d\tau\\
&&=-\int_0^t(\mathbf{v}\cdot\nabla J_{\epsilon}\mathbf{u}, v)d\tau-(\mathbf{v_0},J_{\epsilon}\mathbf{u_0})\notag.
\een

Then by adding up (\ref{sta1}) and (\ref{sta2}) and applying the equality
$$\int_0^{t}(\mathbf{u},\p_tJ_{\epsilon}\mathbf{v})d\tau=-\int_0^{t}(\mathbf{v},\p_tJ_{\epsilon}\mathbf{u})d\tau.$$
one has
\beno
&&-(\mathbf{u},J_{\epsilon}\mathbf{v})(t)-(\mathbf{v},J_{\epsilon}\mathbf{u})(t)-\nu\int_0^t[(\nabla_h \mathbf{u},\nabla_h J_{\epsilon}\mathbf{v})+(\nabla_h \mathbf{v},\nabla_h J_{\epsilon}\mathbf{u})]d\tau\\
&&=-\int_0^t[(\mathbf{u}\cdot\nabla J_{\epsilon}\mathbf{v}, \mathbf{u})+(\mathbf{v}\cdot\nabla J_{\epsilon}\mathbf{u}, v)]d\tau-(\mathbf{u_0,}J_{\epsilon}\mathbf{v_0})-(\mathbf{v_0},J_{\epsilon}\mathbf{u_0})\notag,
\eeno
which implies, after letting $\epsilon\rightarrow0$, that
\ben\label{sta4}
&&-2(\mathbf{u},\mathbf{v})(T)-2\nu\int_0^t(\nabla_h \mathbf{u},\nabla_h \mathbf{v})d\tau\notag\\
&=&-2\int_0^t((\mathbf{v}-\mathbf{u})\cdot\nabla (\mathbf{v}-\mathbf{u}), \mathbf{u})d\tau-2(\mathbf{u_0},\mathbf{v_0}).
\een

Subsequently, by summing up \eqref{sta4} and energy inequalities \eqref{weak1} for $\mathbf{u}$ and $\mathbf{v}$, there holds that
\ben\label{sta5}
\|\mathbf{w}\|_{{L^2({\mathcal{D}})}}^2(t)+\nu\int_0^t\|\nabla_{h}\mathbf{w}\|_{{L^2({\mathcal{D}})}}^2d\tau \leq\|\mathbf{w_0}\|_{{L^2({\mathcal{D}})}}^2+2 \int_0^t(\mathbf{w}\cdot\nabla \mathbf{w}, \mathbf{u})d\tau,
\een
where $\mathbf{w}=\mathbf{v}-\mathbf{u}$. Until now, it suffices to deal with the nonlinear term in (\ref{sta5}). In general, it should be a hard term. Initially, by integrating by parts, using the two-dimensional Ladyzhenskaya inequality in $B_1$ and Young inequality, one has
\ben\label{sta6}
\int_0^t(\mathbf{w}\cdot\nabla \mathbf{w}, \mathbf{u})d\tau\notag
&=&- \int_0^t(\mathbf{w}\cdot\nabla \mathbf{u}, \mathbf{w})d\tau\notag\\
&\leq&\int_0^t\int_0^{2\pi}\|\mathbf{w}\|_{L^4(B_1)}^2\|\nabla \mathbf{u}\|_{L^2(B_1)}d{x_3}d\tau\\
&\leq&\sqrt{2}\int_0^t\int_0^{2\pi}\|\mathbf{w}\|_{L^2(B_1)}\|\nabla_h \mathbf{w}\|_{L^2(B_1)}\|\nabla \mathbf{u}\|_{L^2(B_1)}d{x_3}d\tau\notag\\
&\leq&\f{\nu}2\int_0^t\|\nabla_h \mathbf{w}\|_{L^2(\mathcal{D})}^2d\tau+\f1{\nu}\int_0^t\int_0^{2\pi}\|\mathbf{w}\|_{L^2(B_1)}^2\|\nabla \mathbf{u}\|_{L^2(B_1)}^2d{x_3}d\tau.\notag
\een

Hence, by using Lemma \ref{twointegral}, we can update \eqref{sta6} as
\beno
\int_0^t(\mathbf{w}\cdot\nabla \mathbf{w}, \mathbf{u})d\tau\leq\f{\nu}2\int_0^t\|\nabla_h \mathbf{w}\|_{L^2(\mathcal{D})}^2d\tau+\f{ c^2_*}{\nu}\int_0^t\|\mathbf{w}\|_{L^2(\mathcal{D})}^2\|\mathbf{u}\|_{H^1(\mathcal{D})}^2d\tau,
\eeno
which together with Gronwall's inequality and Corollary \ref{lem1} yield that
\beno
\|\mathbf{w}\|_{{L^2({\mathcal{D}})}}^2(t)&\leq&\|\mathbf{w_0}\|_{{L^2({\mathcal{D}})}}^2\,{\rm exp} \big[\f{c_*^2}{\nu}\int_0^t\| \mathbf{u}\|_{H^1(\mathcal{D})}^2d\tau\big]\\
&\leq&\|\mathbf{w_0}\|_{{L^2({\mathcal{D}})}}^2\,{\rm exp} \big[\f{4c_*^2+2c_*^2c_0}{{\nu^2}}\|\mathbf{u_0}\|_{L^2(\mathcal{D})}^2\big].
\eeno
\endproof\vskip 0.5cm

\subsection{Proof of Theorem \ref{gs}}\hspace*{\parindent}\vskip 0.5cm

With the help of Theorem \ref{gw}, to prove Theorem \ref{gs}, it suffices to verify that
$\mathbf{u}\in L^\infty(0,t;H^1(\mathcal{D}))$ and $\nabla_{h}\mathbf{u}\in L^2(0,t;H^1(\mathcal{D}))$
hold for any $t>0$. By applying Lemma \ref{lem3}, one can achieve the conclusion.
\endproof

\section*{Acknowledgments}
Niu is supported by National Natural Sciences Foundation of China (No. 11471220), the Beijing Natural Science Foundation grants (No. 1142004) and Beijing Municipal Commission of Education grants (No. KM201610028001).
\par


\begin{thebibliography}{10}
\bibitem{AnBr} M. Andersen and M. Br$\o$ns, {\it Topology of helical fluid flow}, European J.
Appl. Mech. 25 (2014), 275-396.

\bibitem{BMNHT} C. Bardos, M. C. Lopes Filho, D. Niu, H. Nussenzveig Lopes, E. S. Titi, {\it Stability of two-dimensional viscous incompressible flows under three-dimensional perturbations and inviscid symmetry breaking}, SIAM J. Math. Anal. 45 (2013), no. 3, 1871-1885.

\bibitem{BLL} A. Bronzi, M. Lopes and H. Lopes Nuzzenveig. {\it Global existence of a weak solution
of the incompressible Euler equations with helical symmetry and $L^p$ vorticity},
Indiana Univ. Math. J. 64, (2015), no. 1, 309-341.

\bibitem{CJT1} C. Cao, J. Li, E. S. Titi, {\it Global well-posedness of strong solutions to the 3D primitive equations with horizontal eddy diffusivity}, J. Differential Equations 257 (2014), no. 11, 4108-4132.

\bibitem{CaoWu1}  C. Cao, J. Wu, {\it Global regularity for the 2D MHD equations with mixed partial dissipation and magnetic diffusion}, Adv. Math. 226 (2011), no. 2, 1803-1822.

\bibitem{CaoWu3} C. Cao, J. Wu, {\it Global regularity for the two-dimensional anistropic Boussinesq equations with vertical dissipation}, Arch. Ration. Mech. Anal. 208 (2013), no. 3, 985-1004.

\bibitem{Danchin} P. Danchin, M. Paicu, {\it Global existence results for the anistropic Boussinesq system in dimension two}, Math. Models Methods Appl. Sci. 21 (2011), no. 3, 421-457.

\bibitem{Du} A. Dutrifoy, {\it Existence globale en temps de solutions h{\'e}lico{\"i}dales des {\'e}quations d'Euler}, C. R. Acad. Sci. Paris S{\'e}r I Math., 329(7):653-656, 1999.

\bibitem{ET} B. Ettinger, E. S. Titi, {\it Global existence and uniqueness of weak solutions of 3D Euler equations with helical symmetric in the absence of vorticity stretching}, SIAM
J. Math. Anal. 41(1) (2009), 269-296.

\bibitem{JLN} Quansen Jiu, Jun Li, Dongjuan Niu, {\it Global Existence of Weak Solutions to the Three-dimensional Euler
Equations with Helical Symmetry}, Preprint, 2016.

\bibitem{JT1} J. Li, E. S. Titi, {\it Global Well-Posedness of the 2D Boussinesq Equations with Vertical Dissipation}, Arch. Ration. Mech. Anal. 220 (2016), no. 3, 983-1001.

\bibitem{MLLNT} M. Lopes, A. Mazzucato, D. Niu, H. Lopes Nussenzveig, E. Titi, {\it Planar limits of
incompressible viscous flows with helical symmetry.} Journal of Dynamics and Diffential
Equations 26, (2014), 843-869.


\bibitem{MTL}A. Mahalov, E. S. Titi, S. Leibovich, {\it Invariant helical subspaces for the Navier-Stokes equations}, Arch.
Ration. Mech. Anal. 112 (1990), no. 3, 193-222.

\bibitem{Majda} A. Majda and A. Bertozzi,  {\it Vorticity and Incompressible Flow}, Cambridge University Press, Cambridge, UK, 2002.


\bibitem{MZ1}  C. Miao, X. Zheng, {\it Global well-posedness for axisymmetric Boussinesq system with horizontal viscosity}, J. Math. Pures Appl. (9) 101 (2014), no. 6, 842-872.

\bibitem{MZ2}  C. Miao, X. Zheng, {\it On the global well-posedness for the Boussinesq system with horizontal dissipation}, Comm. Math. Phys. 321 (2013), no. 1, 33-67.

\bibitem{RT} R. Temam, {\it Navier-Stokes Equations: Theory and Numerical Analysis}, North-Holland, Amsterdam, 1984.

\bibitem{WZ} G. Wu, X. Zheng, {\it Global well-posedness for the two-dimensional nonlinear Boussinesq equations with vertical dissipation}, J. Differential Equations 255 (2013), no. 9, 2891-2926.


\end{thebibliography}
\end{document}